\documentclass[a4paper, final, 12pt]{article}
\usepackage{eurosym}
\usepackage{amsmath, amssymb, latexsym, amscd, amsthm,amsfonts,amstext}
\usepackage[mathscr]{eucal}
\usepackage{graphicx}
\usepackage{subfig}
\usepackage{float}
\usepackage{color}
\usepackage{hyperref}
\usepackage[utf8]{inputenc}
\usepackage[english]{babel}
\usepackage{authblk}

\numberwithin{equation}{section}
\begin{document}

\title{Convexification With Viscosity Term for an Inverse Problem of Tikhonov%
}
\author{ Michael V. Klibanov \and Department of Mathematics and Statistics
\and University of North Carolina at Charlotte, Charlotte, NC 28223, USA
\and Email: mklibanv@charlotte.edu}
\date{}
\maketitle

\begin{abstract}
In 1965 A.N. Tikhonov, the founder of the theory of Ill-Posed and Inverse
Problems, has posed an coefficient inverse problem of the recovery of the
unknown electric conductivity coefficient from measurements of the back
reflected electrical signal. In the geophysical application targeted by
Tikhonov, this coefficient depends only on the depth and characterizes the
electrical conductivity of the ground. The goal of this paper is to
construct for this problem a version of the globally convergent
convexification numerical method for this problem. In this version, the
viscosity term is used in the convexification method. A Carleman estimate
allows to prove global convergence of this method.
\end{abstract}



\section{Introduction}

\label{sec:1}

The goal of this paper is to construct a globally convergent numerical
method for a 1d Coefficient Inverse Problem (CIP) and carry out its
convergence analysis. This CIP has an application in geophysics, as was
first noticed in 1965 by A.N. Tikhonov \cite{Tikh}, the founder of the
theory of Ill-Posed and Inverse Problems. This is the problem of the
recovery of the spatial dependence of the electrical conductivity of the
medium under the ground surface from measurements of the backscattering
electric field on that surface. In \cite{Tikh}, so as in this paper, the
electric conductivity coefficient $\sigma \left( z\right) $ is assumed to be
dependent only on the depth $z>0$. This is a reasonable assumption in many
geophysical applications. Tikhonov has proven uniqueness theorem of this CIP
under the assumption that $\sigma \left( z\right) $ is a piecewise analytic
function \cite{Tikh}.

CIPs are both ill-posed and nonlinear. Conventional numerical methods for
CIPs are based on the minimization of least squares mismatch functionals: we
refer to, e.g. \cite{B1,B2,B3,Chavent,Giorgi,Gonch1,Gonch2,Grote,LB,Riz} and
references cited therein for these methods. However, due to the
ill-posedness and the nonlinearity of CIPs, these functionals are non
convex. The non convexity, in turn causes the well known phenomenon of
multiple local minima and ravines of these functionals, see, e.g. \cite%
{Scales} for a numerical example of multiple local minima for a 1d CIP for a
hyperbolic PDE. Gradient-like methods of the minimization of the least
squares cost functionals for CIPs can stop at any point of a local minimum.
Therefore, convergence of such a method to the correct solution of a CIP can
be rigorously guaranteed only if its starting point is located in a
sufficiently small neighborhood of this solution. In other words, a good
first guess about the solution needs to be available. However, it is unclear
in many applications how to obtain such a guess.

\textbf{Definition}. \emph{We call a numerical method for a CIP globally
convergent, if a theorem is proven, which claims its convergence to the true
solution of that problem without any advanced knowledge of a small
neighborhood of this solution.}

The work \cite{KI} is the first one where the convexification concept was
introduced in the field of CIP. The goal of the convexification is to avoid
the phenomenon of local minima. The convexification is a general concept of
construction of globally convergent numerical methods for CIPs. Each CIP\
requires it own version of the convexification method with its own
convergence analysis. Currently various versions of the convexification
method are developed for CIPs for three main types of PDEs of the second
order: elliptic, parabolic and hyperbolic ones as well as for the radiative
transport equation, see, e.g. \cite{KL,SAR,Khyp,Ktransp} and references
cited therein. In particular, the book \cite{KL} contains results obtained
prior its publication in 2021. The convexification works only for formally
determined CIPs, in which the number $m$ of independent variables in the
input data equals the number $n$ of independent variables in the unknown
coefficient, $m=n$. In our case $m=n=1.$

We introduce the viscosity term in our numerical scheme and develop the
convexification method for this case. For the first time, the  viscosity
term was introduced in the convexification method in \cite{KlibHJ} to
numerically solve the Hamilton-Jacobi equation. We also refer to \cite{NN}
in this regard. In \cite{Ktransp} the viscosity term was used for the
convexification method for a CIP for the radiative transport equation.

The convexification method consists of two steps. On the first step, called
\textquotedblleft transformation", the CIP is transformed to a boundary
value problem (BVP) for such a PDE (or a system of PDEs), which does not
contain the unknown coefficient. Both Dirichlet and Neumann boundary
conditions are given for this BVP. This BVP\ is solved on the second step.
On this step a weighted Tikhonov-like functional is constructed. The weight
is the Carleman Weight Function (CWF). This is the function, which is
involved as the weight in the Carleman estimate for the corresponding PDE
operator, see, e.g. books \cite{KL,LRS} for Carleman estimates. The key
theorem of our convergence analysis claims that this functional is strongly
convex on a convex bounded set of an appropriate Hilbert space. Since a
smallness condition is not imposed on the diameter $d>0$ of this set, then
this is the global strong convexity. The final theorem of the convergence
analysis claims the convergence to the true solution of the CIP of the
gradient descent method of the minimization of that functional, as long as
its starting point is located on the subset of that set, whose diameter is $%
d/6$ and the noise level in the data tends to zero. Again, since a smallness
condition is not imposed on $d$, then this is the global convergence, as
defined above. Explicit convergence rates are also derived.

In section 2 we state both forward and inverse problems. In section 3 we
describe our transformation procedure and construct the above mentioned
functional. In section 4 construct our convexified functional. In section 5
we formulate theorems of our convergence analysis. These theorems are proven
in section 6 and 7.

\section{Forward and Inverse Problems}

\label{sec:2}

Let the number $T>0.$ For $\alpha \in \left( 0,1\right) $ let $C^{\alpha
}\left( \mathbb{R}\right) $ and $C^{2+\alpha ,1+\alpha /2}\left( \mathbb{R}%
\times \left[ 0,T\right] \right) $ be H\"{o}lder spaces \cite{Lad,LU}. For $%
z\in \mathbb{R}$ let $\left\{ z<0\right\} $ be the air and $\left\{
z>0\right\} $ be the ground. Let the electric conductivity coefficient $%
\sigma \left( z\right) $ has the following properties:%
\begin{equation}
\left. 
\begin{array}{c}
\sigma \left( z\right) \in C^{\alpha }\left( \mathbb{R}\right) , \\ 
\sigma \left( z\right) \geq 1,\forall z\in \mathbb{R}, \\ 
\sigma \left( z\right) =1\text{ for }z\in \left\{ z<0\right\} \cup \left\{
z>Z\right\} ,%
\end{array}%
\right.  \label{2.1}
\end{equation}%
where the number $Z>0.$ Let $t>0$ be time and let the function $u\left(
z,t\right) $ be the voltage in a layered medium.

Then $u\left( z,t\right) $ solves the following problem \cite{Tikh} 
\begin{equation}
\left. 
\begin{array}{c}
\sigma \left( z\right) u_{t}=u_{zz},z\in \mathbb{R},t>0, \\ 
u\left( z,0\right) =\delta \left( z\right) .%
\end{array}%
\right.  \label{2.2}
\end{equation}%
We refer to \cite[\S 11-\S 13 of Chapter 4]{Lad} for such a problem for a
general parabolic operator of the second order. In particular, existence and
uniqueness of the solution 
\begin{equation*}
u\in C^{2+\alpha ,1+\alpha /2}\left( \left( \mathbb{R}\times \left[ 0,T%
\right] \right) \diagdown \left\{ \left\vert z\right\vert <\varepsilon ,t\in
\left( 0,\varepsilon \right) \right\} \right) ,\text{ }\forall T>0.
\end{equation*}%
of problem (\ref{2.2}) follows from this reference, where $\varepsilon \in
\left( 0,1\right) $ is an arbitrary number.

\textbf{Coefficient Inverse Problem 1 (CIP1). }\emph{Assume that the
function }$\sigma \left( z\right) $ \emph{satisfies conditions (\ref{2.1})
and the following function }$f\left( t\right) $\emph{\ is known for all }$%
t>0:$%
\begin{equation}
u_{z}\left( 0,t\right) =f\left( t\right) .  \label{2.3}
\end{equation}%
\emph{Find the function }$\sigma \left( z\right) $\emph{\ for }$z\in \left(
0,Z\right) .$

Since we work below only with the Laplace transform of the function $u\left(
z,t\right) ,$ then it is convenient to reformulate CIP1 for the case of the
Laplace transform. For $k>0,$ let 
\begin{equation}
v\left( z,k\right) =\int\limits_{0}^{\infty }u\left( z,t\right) e^{-kt}dt
\label{2.4}
\end{equation}%
be the Laplace transform of the function $u\left( z,t\right) .$ Let $g\left(
k\right) $ be the Laplace transform (\ref{2.4}) of functions $f_{0}\left(
t\right) $ and $f_{1}\left( t\right) $ in (\ref{2.3}) respectively. Then
CIP1 is transformed in CIP2, and we work below only with CIP2.

\textbf{Coefficient Inverse Problem 2 (CIP2). }\emph{Assume that the
function }$\sigma \left( z\right) $ \emph{satisfies conditions (\ref{2.1})
and the following function }$g\left( k\right) $\emph{\ is known}%
\begin{equation}
v_{z}\left( 0,k\right) =g\left( k\right) ,\text{ }\forall k\geq k_{\min }>0,
\label{2.5}
\end{equation}%
\emph{where }$k_{\min }>0$\emph{\ is a certain number chosen below. Find the
function }$\sigma \left( z\right) $\emph{\ for }$z\in \left( 0,Z\right) .$

\section{Transformation}

\label{sec:3}

It follows from (\ref{2.3}) and (\ref{2.4}) that 
\begin{equation}
v_{zz}-kv-k\left( \sigma \left( z\right) -1\right) v=-\delta \left( z\right)
,z\in \mathbb{R}.  \label{3.1}
\end{equation}%
It follows from the third line of (\ref{2.1}) that we have $v_{zz}-kv=0$ for 
$z>M$ and $z<0.$ Hence, 
\begin{equation*}
v\left( z,k\right) =C_{1}e^{-\sqrt{k}z}+C_{2}e^{\sqrt{k}z},\text{ }z>M,
\end{equation*}%
\begin{equation*}
v\left( z,k\right) =C_{3}e^{\sqrt{k}z},\text{ }z<0.
\end{equation*}%
where the numbers $C_{1},C_{2},C_{3}$ are independent on $z$. Hence, to have
the function $v\left( z,k\right) $ bounded, we set 
\begin{equation}
\left. 
\begin{array}{c}
v\left( z,k\right) =A_{1}e^{-\sqrt{k}z},\text{ }z>M, \\ 
v\left( z,k\right) =A_{2}e^{\sqrt{k}z},\text{ }z<0,%
\end{array}%
\right.  \label{3.2}
\end{equation}%
where the numbers $A_{1},A_{2}$ are independent on $z$.

The fundamental solution $u^{0}\left( z,k\right) $ of equation 
\begin{equation}
u_{zz}^{0}-ku^{0}=-\delta \left( z\right)  \label{3.20}
\end{equation}
is%
\begin{equation}
u^{0}\left( z,k\right) =\frac{\exp \left( -\sqrt{k}\left\vert z\right\vert
\right) }{2\sqrt{k}}.  \label{3.3}
\end{equation}%
By (\ref{3.2}) and (\ref{3.3}) 
\begin{equation}
v\left( 0,k\right) =u^{0}\left( 0,k\right) =\frac{1}{2\sqrt{k}}.  \label{3.4}
\end{equation}%
Replace the function $v\left( z,k\right) $ with the function $w\left(
z,k\right) $, where 
\begin{equation}
v\left( z,k\right) =w\left( z,k\right) u^{0}\left( z,k\right) .  \label{3.5}
\end{equation}%
Using (\ref{3.1}) and (\ref{3.20})-(\ref{3.5}), we obtain%
\begin{equation}
w_{zz}-2\sqrt{k}w_{z}-k\left( \sigma \left( z\right) -1\right) w=0\text{ for 
}z>0.  \label{3.6}
\end{equation}%
Since $u>0$ as the fundamental solution of the parabolic equation in the
first line of (\ref{2.2}) \cite{Lad}, then (\ref{3.3}) and (\ref{3.5}) imply
that $w>0$ as well. Hence, we can consider another change of variables,%
\begin{equation}
p\left( z,k\right) =\frac{1}{k}\ln w\left( z,k\right) .  \label{3.7}
\end{equation}%
Then (\ref{2.5}), (\ref{3.2}) and (\ref{3.3})-(\ref{3.6}) imply%
\begin{equation}
p_{zz}+kp_{z}^{2}-2\sqrt{k}p_{z}=\left( \sigma \left( z\right) -1\right) ,%
\text{ }z\in \left( 0,Z\right) ,  \label{3.8}
\end{equation}%
\begin{equation}
p\left( 0,k\right) =1,  \label{3.9}
\end{equation}%
\begin{equation}
p_{z}\left( 0,k\right) =2\sqrt{k}g\left( k\right) +\frac{4}{k^{3/2}},
\label{3.10}
\end{equation}%
\begin{equation}
p_{z}\left( Z,k\right) =0.  \label{3.11}
\end{equation}%
Differentiate (\ref{3.8}) with respect to $k$ and denote%
\begin{equation}
q\left( z,k\right) =\frac{\partial }{\partial k}p\left( z,k\right) .
\label{3.12}
\end{equation}%
Also, use (\ref{3.9})-(\ref{3.11}). We obtain%
\begin{equation}
q_{zz}+2kq_{z}p_{z}+p_{z}^{2}-2\sqrt{k}q_{z}-\frac{p_{z}}{\sqrt{k}}=0,
\label{3.13}
\end{equation}%
\begin{equation}
q\left( 0,k\right) =0,  \label{3.14}
\end{equation}%
\begin{equation}
q_{z}\left( 0,k\right) =2\sqrt{k}g^{\prime }\left( k\right) +\frac{g\left(
k\right) }{\sqrt{k}}-\frac{6}{k^{5/2}},  \label{3.15}
\end{equation}%
\begin{equation}
q_{z}\left( Z,k\right) =0.  \label{3.16}
\end{equation}

We need now to solve problem (\ref{3.13})-(\ref{3.16}). However, equation (%
\ref{3.13}) contains two unknown functions $p$ and $q.$ Since the initial
condition at any value of $k$ is unknown for the function $p,$ then $p$
cannot be expressed via $q$ using (\ref{3.12}). Hence, we introduce the
viscosity term in equation (\ref{3.13}). More precisely, we perturb this
equation by the viscosity term $-\varepsilon p_{zz},$ where the small
parameter $\varepsilon >0$ needs to be found numerically. We obtain 
\begin{equation}
-\varepsilon p_{zz}+q_{zz}+2kq_{z}p_{z}+p_{z}^{2}-2\sqrt{k}q_{z}-\frac{p_{z}%
}{\sqrt{k}}=0,  \label{3.17}
\end{equation}%
We cannot prove convergence of the procedure described below as $\varepsilon
\rightarrow 0.$ It is well known that such a proof is a very non-trivial one
for any PDE. Hence, we do not address this question in the current paper, so
as in two previous publications of the first author with coauthors about the
convexification method with the viscosity term \cite{KlibHJ,Ktransp}.

Introduce a new function $r\left( z,k,\varepsilon \right) ,$ 
\begin{equation}
r\left( z,k,\varepsilon \right) =q-\varepsilon p.  \label{3.18}
\end{equation}%
Hence,%
\begin{equation}
p=\frac{q-r}{\varepsilon }.  \label{3.19}
\end{equation}%
Substituting (\ref{3.18}) and (\ref{3.19}) in equations (\ref{3.13}) and (%
\ref{3.17}), we obtain 
\begin{equation}
L_{1}\left( q,r\right) =q_{zz}+2\frac{k}{\varepsilon }q_{z}\left(
q_{z}-r_{z}\right) +\frac{1}{\varepsilon ^{2}}\left( q_{z}-r_{z}\right) ^{2}-
\label{3.200}
\end{equation}%
\begin{equation*}
-2\sqrt{k}q_{z}-\frac{\left( q_{z}-r_{z}\right) }{\varepsilon \sqrt{k}}=0,
\end{equation*}%
\begin{equation}
L_{2}\left( q,r\right) =r_{zz}+2\frac{k}{\varepsilon }q_{z}\left(
q_{z}-r_{z}\right) +\frac{1}{\varepsilon ^{2}}\left( q_{z}-r_{z}\right) ^{2}-
\label{3.21}
\end{equation}%
\begin{equation*}
-2\sqrt{k}q_{z}-\frac{\left( q_{z}-r_{z}\right) }{\varepsilon \sqrt{k}}=0.
\end{equation*}%
Next, (\ref{3.9})-(\ref{3.11}), (\ref{3.14})-(\ref{3.16}) and (\ref{3.18})
imply the following boundary conditions for functions $q\left( z,k\right) $
and $r\left( z,k\right) :$ 
\begin{equation}
q\left( 0,k\right) =0,  \label{3.22}
\end{equation}%
\begin{equation}
q_{z}\left( 0,k\right) =2\sqrt{k}g^{\prime }\left( k\right) +\frac{g\left(
k\right) }{\sqrt{k}}-\frac{6}{k^{5/2}},  \label{3.23}
\end{equation}%
\begin{equation}
q_{z}\left( Z,k\right) =0,  \label{3.24}
\end{equation}%
\begin{equation}
r\left( 0,k\right) =-\varepsilon ,  \label{3.25}
\end{equation}%
\begin{equation}
r_{z}\left( 0,k\right) =2\sqrt{k}\left( g^{\prime }\left( k\right)
-\varepsilon g\left( k\right) \right) +\frac{g\left( k\right) }{\sqrt{k}}-%
\frac{6}{k^{5/2}}-\frac{4\varepsilon }{k^{3/2}},  \label{3.26}
\end{equation}%
\begin{equation}
r_{z}\left( Z,k\right) =0.  \label{3.27}
\end{equation}

Thus, we focus below on the numerical solution of the boundary value problem
(\ref{3.200})-(\ref{3.27}). We explain in the next section how to get an
approximation for the unknown coefficient $\sigma \left( z\right) $ using
the solution of problem (\ref{3.200})-(\ref{3.27}).

\section{Convexification}

\label{sec:4}

We construct in this section a globally convergent numerical method for
problem (\ref{3.200})-(\ref{3.27}). Fix an arbitrary number $R>0$. The small
parameter $\varepsilon >0$ is fixed below. Let $k\in \left[ k_{\min
},k_{\max }\right] ,$ where $k_{\min },k_{\max }$ are two numbers, which we
should choose numerically. We solve problem (\ref{3.200})-(\ref{3.27}) on
the following set $\overline{B\left( R\right) }:$%
\begin{equation}
\overline{B\left( R\right) }=\left\{ 
\begin{array}{c}
\left( q\left( z,k\right) ,r\left( z,k\right) \right) \in H^{2}\left(
0,Z\right) \times H^{2}\left( 0,Z\right) ,\text{ }\forall k\in \left[
k_{\min },k_{\max }\right] , \\ 
\left\Vert q\left( z,k\right) \right\Vert _{H^{2}\left( 0,Z\right)
}+\left\Vert r\left( z,k\right) \right\Vert _{H^{2}\left( 0,Z\right) }\leq R,%
\text{ }\forall k\in \left[ k_{\min },k_{\max }\right] , \\ 
\text{functions }q\text{ and }r\text{ satisfy } \\ 
\text{boundary conditions (\ref{3.22})-(\ref{3.27})}%
\end{array}%
\right\} .  \label{4.1}
\end{equation}%
The sets like $\overline{B\left( R\right) }$ are called \textquotedblleft
correctness sets" in the regularization theory \cite{TA,T}.

Consider the following Carleman Weight Function (CWF):%
\begin{equation}
\varphi _{\lambda }\left( z\right) =e^{-2\lambda z},  \label{4.2}
\end{equation}%
where $\lambda \geq 1$ is a parameter. We now formulate a Carleman estimate
for the operator $d^{2}/dz^{2}.$ First, we define the subspace $%
H_{0}^{2}\left( 0,Z\right) $ of the space $H^{2}\left( 0,Z\right) $ as: 
\begin{equation*}
H_{0}^{2}\left( 0,Z\right) =\left\{ u\left( z\right) \in H^{2}\left(
0,Z\right) :u\left( 0\right) =u^{\prime }\left( 0\right) =0\right\} .
\end{equation*}

\textbf{Theorem 4.1} (Carleman estimate \cite{K2017}). \emph{There exists a
sufficiently large number }$\lambda _{0}=\lambda _{0}\left( Z\right) \geq 1$%
\emph{\ and a number }$C_{0}=C_{0}\left( Z\right) >0,$\emph{\ both numbers
depend only on }$Z$\emph{, such that the following Carleman estimate holds
for the operator }$d^{2}/dz^{2}$%
\begin{equation}
\int\limits_{0}^{Z}u_{zz}^{2}\varphi _{\lambda }dz\geq
C_{0}\int\limits_{0}^{Z}u_{zz}^{2}\varphi _{\lambda }dz+C_{0}\lambda
\int\limits_{0}^{Z}\left( u_{z}^{2}+\lambda ^{2}u^{2}\right) \varphi
_{\lambda }dz,  \label{4.3}
\end{equation}%
\begin{equation*}
\forall u\in H_{0}^{2}\left( 0,Z\right) ,\text{ }\forall \lambda \geq
\lambda _{0}.
\end{equation*}

Let $L_{1}\left( q,r\right) $ and $L_{2}\left( q,r\right) $ be the nonlinear
partial differential operators defined in (\ref{3.200}) and (\ref{3.21}).
Consider the weighted Tikhonov-like functional, where the weight function is
the CWF of (\ref{4.2}),%
\begin{equation}
\left. J_{\lambda }\left( q,r\right) \left( k\right) =\int\limits_{0}^{Z} 
\left[ \left( L_{1}\left( q,r\right) \left( z,k\right) \right) ^{2}+\left(
L_{2}\left( q,r\right) \left( z,k\right) \right) ^{2}\right] \varphi
_{\lambda }\left( z\right) dz,\right.  \label{4.4}
\end{equation}%
where $\alpha \in \left( 0,1\right) $ is the regularization parameter. The
functional $J_{\lambda ,\alpha }\left( q,r\right) $ depends on $k\in \left[
k_{\min },k_{\max }\right] $ since functions $q,r$ depend on $k$. Here $%
\left[ k_{\min },k_{\max }\right] $ is a certain interval, which should be
chosen computationally.

We solve the following Minimization Problem:

\textbf{Minimization Problem.}\emph{\ For each }$k\in \left[ k_{\min
},k_{\max }\right] ,$\emph{\ minimize functional (\ref{4.4}) on the set }$%
\overline{B\left( R\right) }$\emph{\ defined in (\ref{4.1}).}

The solution of the problem depends on $k,\lambda $ as on parameters. We
will fix an optimal $\lambda \geq \lambda _{0}$ but will vary $k\in \left[
k_{\min },k_{\max }\right] .$ Suppose that a minimizer $\left( q_{\min
,\lambda }\left( z,k\right) ,r_{\min ,\lambda }\left( z,k\right) \right) $
of functional (\ref{4.4}) is found. Then, using (\ref{3.19}), we set%
\begin{equation}
p_{\min ,\lambda }\left( z,k\right) =\frac{q_{\min ,\lambda }\left(
z,k\right) -r_{\min ,\lambda }\left( z,k\right) }{\varepsilon }.  \label{4.5}
\end{equation}%
Substituting (\ref{4.5}) in (\ref{3.8}), we obtain%
\begin{equation}
\sigma _{\min ,\lambda }\left( z,k\right) =\partial _{z}^{2}p_{\min ,\lambda
}\left( z,k\right) +k\left( \partial _{z}p_{\min ,\lambda }\left( z,k\right)
\right) ^{2}-2\sqrt{k}\partial _{z}p_{\min ,\lambda }\left( z,k\right) +1.
\label{4.6}
\end{equation}%
Finally, we define the computational solution of our Coefficient Inverse
Problem 2 with the input data (\ref{2.5}) as the average over the interval $%
\left[ k_{\min },k_{\max }\right] $ of functions $\sigma _{\min ,\lambda
}\left( z,k\right) $ in (\ref{4.6}), 
\begin{equation}
\sigma _{\text{comp,}\lambda }\left( z\right) =\frac{1}{k_{\max }-k_{\min }}%
\int\limits_{k_{\min }}^{k_{\max }}\sigma _{\min ,\lambda }\left( z,k\right)
dk.  \label{4.7}
\end{equation}

\section{Theorems of the Convergence Analysis}

\label{sec:5}

\subsection{The central result}

\label{sec:5.1}

\textbf{Theorem 5.1} (strong convexity, the central result).

\emph{1. For each }$\lambda >0,$\emph{\ for each }$k>0$\emph{, for each }$%
\alpha >0$\emph{\ and for each point }$\left( q,r\right) \in \overline{%
B\left( R\right) }$\emph{\ the functional }$J_{\lambda ,\alpha }\left(
q,r\right) \left( k\right) $\emph{\ has the Fr\'{e}chet derivative }%
\begin{equation}
\left. 
\begin{array}{c}
J_{\lambda }^{\prime }\left( q,r\right) \left( k\right) \in \left(
H_{0}^{2}\left( 0,Z\right) \times H_{0}^{2}\left( 0,Z\right) \right) \cap \\ 
\cap \left\{ \left( u\left( z,k\right) ,v\left( z,k\right) \right)
:u^{\prime }\left( Z,k\right) =v^{\prime }\left( Z,k\right) =0\right\} .%
\end{array}%
\right.  \label{5.1}
\end{equation}%
\emph{This derivative is Lipschitz continuous on the set }$\overline{B\left(
R\right) }$\emph{, i.e. there exists a number }$D=D\left( R,\lambda
,k\right) >0$\emph{\ such that for all }$k\in \left[ k_{\min },k_{\max }%
\right] $ 
\begin{equation}
\left. 
\begin{array}{c}
\left\Vert J_{\lambda }^{\prime }\left( q_{2},r_{2}\right) \left( k\right)
-J_{\lambda ,\alpha }^{\prime }\left( q_{1},r_{1}\right) \left( k\right)
\right\Vert _{H^{2}\left( 0,Z\right) \times H^{2}\left( 0,Z\right) }\leq \\ 
\leq D\left\Vert \left( q_{2},r_{2}\right) \left( k\right) -\left(
q_{1},r_{1}\right) \left( k\right) \right\Vert _{H^{2}\left( 0,Z\right)
\times H^{2}\left( 0,Z\right) },\text{ }\forall \left( q_{2},r_{2}\right)
,\left( q_{1},r_{1}\right) \in \overline{B\left( R\right) }.%
\end{array}%
\right.  \label{5.2}
\end{equation}

\emph{2. There exists a sufficiently large number }$\lambda _{1}=\lambda
_{1}\left( R,Z,k_{\min },k_{\max },\varepsilon \right) \geq \lambda _{0}\geq
1$\emph{\ such that for all }$\lambda \geq \lambda _{1}$\emph{, for all }$%
k\in \left[ k_{\min },k_{\max }\right] $\emph{\ the functional }$J_{\lambda
,\alpha }\left( q,r\right) $\emph{\ is strongly convex on the set }$%
\overline{B\left( R\right) },$\emph{\ i.e. there exists a number }$%
C_{1}=C_{1}\left( R,Z,k_{\min },k_{\max },\varepsilon \right) >0$\emph{\
such that the following strong convexity inequality is valid: }%
\begin{equation}
\left. 
\begin{array}{c}
J_{\lambda }\left( q_{2},r_{2}\right) \left( k\right) -J_{\lambda }\left(
q_{1},r_{1}\right) \left( k\right) -\left[ J_{\lambda }^{\prime }\left(
q_{1},r_{1}\right) \left( k\right) ,\left( q_{2}-q_{1},r_{2}-r_{1}\right)
\left( k\right) \right] \geq \\ 
\geq C_{1}e^{-2\lambda Z}\left\Vert \left( q_{2},r_{2}\right) \left(
k\right) -\left( q_{1},r_{1}\right) \left( k\right) \right\Vert
_{H^{2}\left( 0,Z\right) \times H^{2}\left( 0,Z\right) }^{2}, \\ 
\forall \left( q_{1},r_{1}\right) \left( k\right) ,\left( q_{2},r_{2}\right)
\left( k\right) \in \overline{B\left( R\right) },\text{ }\forall k\in \left[
k_{\min },k_{\max }\right] ,%
\end{array}%
\right.  \label{5.3}
\end{equation}%
\emph{where }$\left[ ,\right] $\emph{\ is the scalar product in }$%
H^{2}\left( 0,Z\right) \times H^{2}\left( 0,Z\right) .$

\emph{3. For each }$\lambda \geq \lambda _{1}$\emph{, for each }$k\in \left[
k_{\min },k_{\max }\right] $\emph{\ there exists unique minimizer }$\left(
q_{\min ,\lambda }\left( z,k\right) ,r_{\min ,\lambda }\left( z,k\right)
\right) \in \overline{B\left( R\right) }$\emph{\ of the functional} $%
J_{\lambda }\left( q,r\right) \left( k\right) $\emph{\ on the set }$%
\overline{B\left( R\right) }$\emph{\ and the following inequality holds:}%
\begin{equation}
\left. 
\begin{array}{c}
\left[ J_{\lambda ,\alpha }^{\prime }\left( q_{\min ,\lambda },r_{\min
,\lambda }\right) ,\left( q-q_{\min ,\lambda },r-r_{\min ,\lambda }\right) %
\right] \left( k\right) \geq 0,\text{ } \\ 
\forall \left( q,r\right) \left( k\right) \in \overline{B\left( R\right) },%
\text{ }\forall k\in \left[ k_{\min },k_{\max }\right] ,%
\end{array}%
\right.  \label{5.4}
\end{equation}%
\emph{Here and below} $C_{1}=C_{1}\left( R,Z,k_{\min },k_{\max },\varepsilon
\right) >0$\emph{\ denotes different numbers depending only on numbers }$R$%
\emph{, }$Z$\emph{, }$k_{\min }$,$k_{\max }$ \emph{and} $\varepsilon .$ 
\emph{The number }$\lambda _{1}$\emph{\ also depends only on these
parameters.}

\subsection{The accuracy of the minimizer}

\label{sec:5.2}

We now want to estimate the accuracy of the minimizer of functional (\ref%
{4.4}). To do this, we should assume first that there exists exact solution $%
\sigma ^{\ast }\left( z\right) $ satisfying conditions (\ref{2.1}) with the
\textquotedblleft ideal", i.e. noiseless data $g^{\ast }\left( k\right) $ in
(\ref{2.5}). This assumption is a natural one in the theory of Ill-Posed
problems \cite{TA,T}. The function $\sigma ^{\ast }\left( z\right) $
generates functions $q^{\ast }\left( z,k\right) $ and $r^{\ast }\left(
z,k\right) ,$ just as above.

Suppose that there exists a vector function 
\begin{equation*}
F\left( z,k\right) =\left( F_{1},F_{2}\right) \left( z,k\right) \in
H^{2}\left( 0,Z\right) \times H^{2}\left( 0,Z\right) ,\text{ }\forall k\in %
\left[ k_{\min },k_{\max }\right]
\end{equation*}%
satisfying boundary conditions (\ref{3.22})-(\ref{3.27}), where $F_{1}$
stands for $q,$ and $F_{2}$ stands for $r$. Suppose also that there exists a
vector function 
\begin{equation*}
F^{\ast }\left( z,k\right) =\left( F_{1}^{\ast },F_{2}^{\ast }\right) \left(
z,k\right) \in H^{2}\left( 0,Z\right) \times H^{2}\left( 0,Z\right) ,\text{ }%
\forall k\in \left[ k_{\min },k_{\max }\right]
\end{equation*}%
satisfying boundary conditions (\ref{3.22})-(\ref{3.27}), in which $g\left(
k\right) $ is replaced with $g^{\ast }\left( k\right) .$ Again, $F_{1}^{\ast
}$ stands for $q^{\ast }$ and $F_{2}^{\ast }$ stands for $r^{\ast }.$ Let $%
\delta \in \left( 0,1\right) $ be a small number characterizing the noise
level in the data (\ref{3.22})-(\ref{3.27}). More precisely, we assume that 
\begin{equation}
\left\Vert F_{1}-F_{1}^{\ast }\right\Vert _{H^{2}\left( 0,Z\right)
}+\left\Vert F_{2}-F_{2}^{\ast }\right\Vert _{H^{2}\left( 0,Z\right)
}<\delta ,\text{ }\forall k\in \left[ k_{\min },k_{\max }\right] .
\label{5.5}
\end{equation}%
We also assume that%
\begin{equation}
\left. 
\begin{array}{c}
\left\Vert F_{1}\left( z,k\right) \right\Vert _{H^{2}\left( 0,Z\right)
}+\left\Vert F_{2}\left( z,k\right) \right\Vert _{H^{2}\left( 0,Z\right)
}\leq R,\text{ }\forall k\in \left[ k_{\min },k_{\max }\right] , \\ 
\left\Vert F_{1}^{\ast }\left( z,k\right) \right\Vert _{H^{2}\left(
0,Z\right) }+\left\Vert F_{2}^{\ast }\left( z,k\right) \right\Vert
_{H^{2}\left( 0,Z\right) }\leq R,\text{ }\forall k\in \left[ k_{\min
},k_{\max }\right] .%
\end{array}%
\right.  \label{5.6}
\end{equation}%
Let $\overline{B^{\ast }\left( R\right) }$ be the following analog of the
set $\overline{B\left( R\right) }$ in (\ref{4.1}) 
\begin{equation}
\overline{B^{\ast }\left( R\right) }=\left\{ 
\begin{array}{c}
\left( q\left( z,k\right) ,r\left( z,k\right) \right) \in H^{2}\left(
0,Z\right) \times H^{2}\left( 0,Z\right) ,\text{ }\forall k\in \left[
k_{\min },k_{\max }\right] , \\ 
\left\Vert q\left( z,k\right) \right\Vert _{H^{2}\left( 0,Z\right)
}+\left\Vert r\left( z,k\right) \right\Vert _{H^{2}\left( 0,Z\right) }\leq R,%
\text{ }\forall k\in \left[ k_{\min },k_{\max }\right] , \\ 
\text{functions }q\text{ and }r\text{ satisfy } \\ 
\text{boundary conditions (\ref{3.22})-(\ref{3.27}), } \\ 
\text{in which }g\left( k\right) \text{ is replaced with }g^{\ast }\left(
k\right) .%
\end{array}%
\right\} .  \label{5.7}
\end{equation}%
We assume that 
\begin{equation}
\left. 
\begin{array}{c}
R-C_{1}\delta >0, \\ 
\left( q^{\ast }\left( z,k\right) ,r^{\ast }\left( z,k\right) \right) \in
B^{\ast }\left( R-C_{1}\delta \right) .%
\end{array}%
\right.  \label{5.8}
\end{equation}

For every vector function $\left( q,r\right) \in \overline{B\left( R\right) }
$ consider the vector function 
\begin{equation}
\left( \widetilde{q},\widetilde{r}\right) =\left( q-F_{1},r-F_{2}\right) .
\label{5.9}
\end{equation}%
Similarly, consider the vector function 
\begin{equation}
\left( \widetilde{q}^{\ast },\widetilde{r}^{\ast }\right) =\left( q^{\ast
}-F_{1}^{\ast },r^{\ast }-F_{2}^{\ast }\right) .  \label{5.10}
\end{equation}%
It follows from (\ref{5.6})-(\ref{5.10}) and triangle inequality that%
\begin{equation}
\left( \widetilde{q},\widetilde{r}\right) ,\left( \widetilde{q}^{\ast },%
\widetilde{r}^{\ast }\right) \in B_{0}\left( 2R\right) ,  \label{5.100}
\end{equation}%
\begin{equation}
\overline{B_{0}\left( 2R\right) }=\left\{ 
\begin{array}{c}
\left( q\left( z,k\right) ,r\left( z,k\right) \right) \in H^{2}\left(
0,Z\right) \times H^{2}\left( 0,Z\right) ,\text{ }\forall k\in \left[
k_{\min },k_{\max }\right] , \\ 
\left\Vert q\left( z,k\right) \right\Vert _{H^{2}\left( 0,Z\right)
}+\left\Vert r\left( z,k\right) \right\Vert _{H^{2}\left( 0,Z\right) }\leq
2R,\text{ }\forall k\in \left[ k_{\min },k_{\max }\right] , \\ 
\text{functions }q\text{ and }r\text{ satisfy } \\ 
\text{zero boundary conditions (\ref{3.22})-(\ref{3.27}).}%
\end{array}%
\right\}  \label{5.11}
\end{equation}
Consider a new functional 
\begin{equation}
\left. 
\begin{array}{c}
I_{\lambda }\left( q,r\right) :\overline{B_{0}\left( 2R\right) }\rightarrow 
\mathbb{R}, \\ 
I_{\lambda }\left( q,r\right) =J_{\lambda }\left( q+F_{1},r+F_{2}\right) .%
\end{array}%
\right.  \label{5.12}
\end{equation}%
By (\ref{5.12}) an obvious analog of Theorem 5.1 is valid for the functional 
$I_{\lambda }\left( q,r\right) .$ However, since $\left(
q+F_{1},r+F_{2}\right) \in \overline{B\left( 3R\right) }$ for $\left(
q,r\right) :\overline{B_{0}\left( 2R\right) },$then we should take here 
\begin{equation}
\lambda _{2}=\lambda _{1}\left( 3R,Z,k_{\min },k_{\max },\varepsilon \right)
.  \label{5.13}
\end{equation}

\textbf{Theorem 5.2} (the accuracy of the minimizer). \emph{Assume that (\ref%
{5.8})-(\ref{5.13}) hold. For any }$\lambda \geq \lambda _{2},$\emph{\ and
for any }$k\in \left[ k_{\min },k_{\max }\right] ,$\emph{\ let }$\left( 
\widetilde{q}_{\min ,\lambda },\widetilde{r}_{\min ,\lambda }\right) \in 
\overline{B_{0}\left( 2R\right) }$\emph{\ be the unique minimizer on the set 
}$\overline{B_{0}\left( 2R\right) }$ \emph{of the functional }$I_{\lambda
}\left( q,r\right) $\emph{\ in (\ref{5.12}), which is guaranteed by Theorem
5.1, i.e.}%
\begin{equation}
I_{\lambda }\left( \widetilde{q}_{\min ,\lambda },\widetilde{r}_{\min
,\lambda }\right) =\min_{\overline{B_{0}\left( 2R\right) }}I_{\lambda
}\left( q,r\right) .  \label{5.130}
\end{equation}%
\emph{\ Define }%
\begin{equation}
\left( \overline{q}_{\min ,\lambda _{2}},\overline{r}_{\min ,\lambda
_{2}}\right) =\left( \widetilde{q}_{\min ,\lambda _{2}}+F_{1},\widetilde{r}%
_{\min ,\lambda _{2}}+F_{2}\right) .  \label{5.14}
\end{equation}%
\emph{Then} 
\begin{equation}
\left( \overline{q}_{\min ,\lambda _{2}},\overline{r}_{\min ,\lambda
_{2}}\right) \in B\left( R\right) .  \label{5.140}
\end{equation}%
\emph{\ Furthermore, the vector function }$\left( \overline{q}_{\min
,\lambda _{2}},\overline{r}_{\min ,\lambda _{2}}\right) $\emph{\ is the
unique minimizer of the functional }$J_{\lambda _{2}}\left( q,r\right) $%
\emph{\ on the set }$\overline{B\left( R\right) },$ \emph{i.e.}%
\begin{equation}
\left( \overline{q}_{\min ,\lambda _{2}},\overline{r}_{\min ,\lambda
_{2}}\right) =\left( q_{\min ,\lambda _{2}},r_{\min ,\lambda _{2}}\right) ,
\label{5.141}
\end{equation}%
\emph{\ and the following accuracy estimate is valid for all }$k\in \left[
k_{\min },k_{\max }\right] $ 
\begin{equation}
\left\Vert q_{\min ,\lambda _{2}}\left( z,k\right) -q^{\ast }\left(
z,k\right) \right\Vert _{H^{2}\left( 0,Z\right) }+\left\Vert r_{\min
,\lambda _{2}}\left( z,k\right) -r^{\ast }\left( z,k\right) \right\Vert
_{H^{2}\left( 0,Z\right) }\leq C_{1}\delta ,  \label{5.15}
\end{equation}%
\begin{equation}
\left\Vert \sigma _{\text{comp},\lambda _{2}}-\sigma ^{\ast }\right\Vert
_{L_{2}\left( 0,Z\right) }\leq C_{1}\delta ,  \label{5.150}
\end{equation}%
\emph{where the function }$\sigma _{\text{comp},\lambda _{2}}\left( z\right) 
$\emph{\ is found via (\ref{4.5})-(\ref{4.7}). }

\subsection{Global convergence of the gradient descent method}

\label{sec:5.3}

\bigskip Similarly with (\ref{5.8}) assume that 
\begin{equation}
\left. 
\begin{array}{c}
R/3-C_{1}\delta >0, \\ 
\left( q^{\ast }\left( z,k\right) ,r^{\ast }\left( z,k\right) \right) \in
B^{\ast }\left( R/3-C_{1}\delta \right) .%
\end{array}%
\right.  \label{5.16}
\end{equation}%
Let the number $\gamma \in \left( 0,1\right) $ and let%
\begin{equation}
\left( q_{0}\left( z,k\right) ,r_{0}\left( z,k\right) \right) \in B\left( 
\frac{R}{3}\right) .  \label{5.17}
\end{equation}%
Let $\lambda _{2}$ be the number defined in (\ref{5.13}). We construct the
gradient descent method as the following sequence:%
\begin{equation}
\left( q_{n},r_{n}\right) =\left( q_{n-1},r_{n-1}\right) -\gamma J_{\lambda
_{2}}^{\prime }\left( q_{n-1},r_{n-1}\right) ,\text{ }n=1,2,...  \label{5.18}
\end{equation}%
Note that since by Theorem 5.1 $J_{\lambda _{2}}^{\prime }\left(
q_{n-1},r_{n-1}\right) $ satisfies (\ref{5.1}), then all terms of sequence (%
\ref{5.18}) have the same boundary conditions (\ref{3.22})-(\ref{3.27}).

\textbf{Theorem 5.3}. \emph{Assume that (\ref{5.16})-(\ref{5.18}). Let }$%
\left( q_{\min ,\lambda _{2}},r_{\min ,\lambda _{2}}\right) $\emph{\ be the
unique minimizer of the functional }$J_{\lambda _{2}}\left( q,r\right) $%
\emph{\ on the set }$\overline{B\left( R\right) },$ \emph{the existence of
which is guaranteed by Theorem 5.1. Then }$\left( q_{\min ,\lambda
_{2}},r_{\min ,\lambda _{2}}\right) \in B\left( R/3\right) .$\emph{\ There
exists a sufficiently small number }$\gamma \in \left( 0,1\right) $\emph{\
and a number }$\theta =\theta \left( \gamma \right) \in \left( 0,1\right) $%
\emph{\ such that all terms of sequence (\ref{5.18}) belong to }$B\left(
R\right) $\emph{\ and the following estimates hold:}%
\begin{equation}
\left. 
\begin{array}{c}
\left\Vert q_{n}-q^{\ast }\right\Vert _{H^{2}\left( 0,Z\right) }+\left\Vert
r_{n}-r^{\ast }\right\Vert _{H^{2}\left( 0,Z\right) }\leq \\ 
\leq C_{1}\delta +\theta ^{n}\left( \left\Vert q_{\min ,\lambda
_{2}}-q_{0}\right\Vert _{H^{2}\left( 0,Z\right) }+\left\Vert r_{\min
,\lambda _{2}}-r_{0}\right\Vert _{H^{2}\left( 0,Z\right) }\right) , \\ 
\text{ }\forall k\in \left[ k_{\min },k_{\max }\right] ,%
\end{array}%
\right.  \label{5.19}
\end{equation}%
\begin{equation}
\left. 
\begin{array}{c}
\left\Vert \sigma _{\text{comp},\lambda _{2},n}-\sigma ^{\ast }\right\Vert
_{L_{2}\left( 0,Z\right) }\leq \\ 
\leq C_{1}\delta +\theta ^{n}\sup_{k\in \left[ k_{\min },k_{\max }\right]
}\left( \left\Vert q_{\min ,\lambda }-q_{0}\right\Vert _{H^{2}\left(
0,Z\right) }+\left\Vert r_{\min ,\lambda }-r_{0}\right\Vert _{H^{2}\left(
0,Z\right) }\right) ,%
\end{array}%
\right.  \label{5.20}
\end{equation}%
\emph{where functions }$\sigma _{\text{comp},\lambda ,n}\left( z\right) $%
\emph{\ are defined as in (\ref{4.7}) with the replacement of the triple }$%
\left( q_{\min ,\lambda },r_{\min ,\lambda },p_{\min ,\lambda }\right) $%
\emph{\ with }$\left( q_{n},r_{n},p_{n}\right) .$

\textbf{Proof.} Assuming that Theorems 5.1 and 5.2 are valid, the proof
follows immediately from Theorem 6 of \cite{SAR}. $\square $

\textbf{Remark 5.1.} \emph{Since smallness conditions are not imposed on }$R$%
\emph{, then Theorem 5.3 claims the global convergence of sequence (\ref%
{5.18}), see Definition in section 1.}

\section{Proof of Theorem 5.1}

\label{sec:6}

Below $C_{2}=C_{2}\left( R,Z,k_{\min },k_{\max },\varepsilon \right) >0$%
\emph{\ }denotes different numbers depending only on numbers $R$, $Z$, $%
k_{\min }$,$k_{\max }$ and $\varepsilon .$ Consider two arbitrary pairs $%
\left( q_{1},r_{1}\right) ,\left( q_{2},r_{2}\right) \in \overline{B\left(
R\right) }.$ Then by (\ref{4.1}), (\ref{5.11}) and triangle inequality 
\begin{equation}
\left( q_{2},r_{2}\right) -\left( q_{1},r_{1}\right) =\left(
h_{1},h_{2}\right) \in \overline{B_{0}\left( 2R\right) }.  \label{6.1}
\end{equation}%
Also, (\ref{4.1}), (\ref{5.11}), (\ref{6.1}) and Sobolev embedding theorem
imply that%
\begin{equation}
\left. 
\begin{array}{c}
\left( q_{1},r_{1}\right) ,\left( q_{2},r_{2}\right) ,\left(
h_{1},h_{2}\right) \in C^{1}\left[ 0,Z\right] \times C^{1}\left[ 0,Z\right] ,
\\ 
\left\Vert \left( q_{1},r_{1}\right) \right\Vert _{C^{1}\left[ 0,Z\right]
\times C^{1}\left[ 0,Z\right] },\left\Vert \left( q_{2},r_{2}\right)
\right\Vert _{C^{1}\left[ 0,Z\right] \times C^{1}\left[ 0,Z\right]
},\left\Vert \left( h_{1},h_{2}\right) \right\Vert _{C^{1}\left[ 0,Z\right]
\times C^{1}\left[ 0,Z\right] }\leq C, \\ 
\forall k\in \left[ k_{\min },k_{\max }\right] ,%
\end{array}%
\right.  \label{6.2}
\end{equation}%
where the number $C=C\left( R,Z,k_{\min },k_{\max }\right) >0$ is a number
depending only on listed parameters. By (\ref{3.200}) and (\ref{6.1}) 
\begin{equation*}
L_{1}\left( q_{2},r_{2}\right) =L_{1}\left( q_{1}+h_{1},r_{1}+h_{2}\right) =
\end{equation*}%
\begin{equation*}
=\partial _{z}^{2}q_{1}+\partial _{z}^{2}h_{1}+2\frac{k}{\varepsilon }\left(
\partial _{z}q_{1}+\partial _{z}h_{1}\right) \left[ \left( \partial
_{z}q_{1}-\partial _{z}r_{1}\right) +\left( \partial _{z}h_{1}-\partial
_{z}h_{2}\right) \right] +
\end{equation*}%
\begin{equation*}
+\frac{1}{\varepsilon ^{2}}\left( \partial _{z}q_{1}-\partial
_{z}r_{1}\right) ^{2}+\frac{2}{\varepsilon ^{2}}\left( \partial
_{z}q_{1}-\partial _{z}r_{1}\right) \left( \partial _{z}h_{1}-\partial
_{z}h_{2}\right) +\frac{2}{\varepsilon ^{2}}\left( \partial
_{z}h_{1}-\partial _{z}h_{2}\right) ^{2}-
\end{equation*}%
\begin{equation*}
-2\sqrt{k}\partial _{z}q_{1}-\frac{\left( \partial _{z}q_{1}-\partial
_{z}r_{1}\right) }{\varepsilon \sqrt{k}}-2\sqrt{k}\partial _{z}h_{1}-\frac{%
\left( \partial _{z}h_{1}-\partial _{z}h_{2}\right) }{\varepsilon \sqrt{k}}.
\end{equation*}%
We now single out the linear, with respect to $\left( h_{1},h_{2}\right) $
part of this expression. We obtain 
\begin{equation*}
L_{1}\left( q_{2},r_{2}\right) =L_{1}\left( q_{1},r_{1}\right) +
\end{equation*}%
\begin{equation*}
+\partial _{z}^{2}h_{1}+2\frac{k}{\varepsilon }\partial _{z}h_{1}+2\frac{k}{%
\varepsilon }\partial _{z}q_{1}\left( \partial _{z}h_{1}-\partial
_{z}h_{2}\right) +\frac{2}{\varepsilon ^{2}}\left( \partial
_{z}q_{1}-\partial _{z}r_{1}\right) \left( \partial _{z}h_{1}-\partial
_{z}h_{2}\right) -
\end{equation*}%
\begin{equation}
-2\sqrt{k}\partial _{z}h_{1}-\frac{\left( \partial _{z}h_{1}-\partial
_{z}h_{2}\right) }{\varepsilon \sqrt{k}}+  \label{6.3}
\end{equation}%
\begin{equation*}
+2\frac{k}{\varepsilon }\partial _{z}h_{1}\left( \partial _{z}h_{1}-\partial
_{z}h_{2}\right) +\frac{2}{\varepsilon ^{2}}\left( \partial
_{z}h_{1}-\partial _{z}h_{2}\right) ^{2}.
\end{equation*}%
Denote 
\begin{equation*}
L_{1,\text{linear}}\left( h_{1},h_{2}\right) =\partial _{z}^{2}h_{1}+2\frac{k%
}{\varepsilon }\partial _{z}h_{1}+2\frac{k}{\varepsilon }\partial
_{z}q_{1}\left( \partial _{z}h_{1}-\partial _{z}h_{2}\right) +
\end{equation*}%
\begin{equation}
+\frac{2}{\varepsilon ^{2}}\left( \partial _{z}q_{1}-\partial
_{z}r_{1}\right) \left( \partial _{z}h_{1}-\partial _{z}h_{2}\right) -2\sqrt{%
k}\partial _{z}h_{1}-\frac{\left( \partial _{z}h_{1}-\partial
_{z}h_{2}\right) }{\varepsilon \sqrt{k}},  \label{6.4}
\end{equation}%
In addition, denote%
\begin{equation}
L_{1,\text{nonlinear}}\left( h_{1},h_{2}\right) =2\frac{k}{\varepsilon }%
\partial _{z}h_{1}\left( \partial _{z}h_{1}-\partial _{z}h_{2}\right) +\frac{%
2}{\varepsilon ^{2}}\left( \partial _{z}h_{1}-\partial _{z}h_{2}\right) ^{2}.
\label{6.5}
\end{equation}%
Using (\ref{6.3})-(\ref{6.5}), we obtain%
\begin{equation*}
\left( L_{1}\left( q_{2},r_{2}\right) \right) ^{2}-\left( L_{1}\left(
q_{2},r_{2}\right) \right) ^{2}=2L_{1}\left( q_{1},r_{1}\right) L_{1,\text{%
linear}}\left( h_{1},h_{2}\right) +
\end{equation*}%
\begin{equation}
+\left( L_{1,\text{linear}}\left( h_{1},h_{2}\right) \right)
^{2}+2L_{1}\left( q_{1},r_{1}\right) L_{1,\text{nonlinear}}\left(
h_{1},h_{2}\right) +  \label{6.6}
\end{equation}%
\begin{equation*}
+2L_{1,\text{linear}}\left( h_{1},h_{2}\right) L_{1,\text{nonlinear}}\left(
h_{1},h_{2}\right) +\left( L_{1,\text{nonlinear}}\left( h_{1},h_{2}\right)
\right) ^{2}.
\end{equation*}%
Next, (\ref{6.2})-(\ref{6.6}) and Cauchy-Schwarz inequality lead to%
\begin{equation*}
\left( L_{1}\left( q_{2},r_{2}\right) \right) ^{2}-\left( L_{1}\left(
q_{1},r_{1}\right) \right) ^{2}-2L_{1}\left( q_{1},r_{1}\right) L_{1,\text{%
linear}}\left( h_{1},h_{2}\right) \geq
\end{equation*}%
\begin{equation}
\geq \frac{1}{2}\left( \partial _{z}^{2}h_{1}\right) ^{2}-C_{1}\left[ \left(
\partial _{z}h_{1}\right) ^{2}+\left( \partial _{z}h_{2}\right) ^{2}\right] .
\label{6.7}
\end{equation}%
Similarly, using (\ref{3.21}), (\ref{6.1}), (\ref{6.2}) and analogs of
formulas (\ref{6.3})-(\ref{6.6}) for the operator $L_{2},$ we obtain the
following analog of (\ref{6.7}):%
\begin{equation*}
\left( L_{2}\left( q_{2},r_{2}\right) \right) ^{2}-\left( L_{2}\left(
q_{1},r_{1}\right) \right) ^{2}-2L_{2}\left( q_{1},r_{1}\right) L_{2,\text{%
linear}}\left( h_{1},h_{2}\right) \geq
\end{equation*}%
\begin{equation}
\geq \frac{1}{2}\left( \partial _{z}^{2}h_{2}\right) ^{2}-C_{1}\left[ \left(
\partial _{z}h_{1}\right) ^{2}+\left( \partial _{z}h_{2}\right) ^{2}\right] .
\label{6.8}
\end{equation}

Hence, using (\ref{4.4}), we obtain 
\begin{equation*}
J_{\lambda }\left( q_{2},r_{2}\right) \left( k\right) -J_{\lambda }\left(
q_{1},r_{1}\right) \left( k\right) =
\end{equation*}%
\begin{equation*}
=2\int\limits_{0}^{Z}\left[ L_{1}\left( q_{1},r_{1}\right) L_{1,\text{linear}%
}\left( h_{1},h_{2}\right) +L_{2}\left( q_{1},r_{1}\right) L_{2,\text{linear}%
}\left( h_{1},h_{2}\right) \right] \varphi _{\lambda }\left( z\right) dz+
\end{equation*}%
\begin{equation}
+\sum\limits_{i=1}^{2}\int\limits_{0}^{Z}\left[ \left( L_{i,\text{linear}%
}\left( h_{1},h_{2}\right) \right) ^{2}+2L_{i}\left( q_{1},r_{1}\right) L_{i,%
\text{nonlinear}}\left( h_{1},h_{2}\right) \right] \varphi _{\lambda }\left(
z\right) dz+  \label{6.9}
\end{equation}%
\begin{equation*}
+\sum\limits_{i=1}^{2}\int\limits_{0}^{Z}\left[ 2L_{i,\text{linear}}\left(
h_{1},h_{2}\right) L_{i,\text{nonlinear}}\left( h_{1},h_{2}\right) +\left(
L_{i,\text{nonlinear}}\left( h_{1},h_{2}\right) \right) ^{2}\right] \varphi
_{\lambda }\left( z\right) dz.
\end{equation*}%
Consider the expression in the second line of (\ref{6.9}),%
\begin{equation*}
\widehat{J}_{\lambda ,q_{1},r_{1}}\left( h_{1},h_{2}\right) \left( k\right) =
\end{equation*}%
\begin{equation}
=2\int\limits_{0}^{Z}\left[ L_{1}\left( q_{1},r_{1}\right) L_{1,\text{linear}%
}\left( h_{1},h_{2}\right) +L_{2}\left( q_{1},r_{1}\right) L_{2,\text{linear}%
}\left( h_{1},h_{2}\right) \right] \varphi _{\lambda }\left( z\right) dz.
\label{6.10}
\end{equation}%
Obviously, 
\begin{equation}
\widehat{J}_{\lambda ,q_{1},r_{1}}\left( h_{1},h_{2}\right) \left( k\right)
:H_{0}^{2}\left( 0,Z\right) \times H_{0}^{2}\left( 0,Z\right) \rightarrow 
\mathbb{R}  \label{6.11}
\end{equation}%
is a bounded linear functional. Hence, by Riesz theorem there exists unique
vector function $\widetilde{J}_{\lambda ,q_{1},r_{1}}\in H_{0}^{2}\left(
0,Z\right) \times H_{0}^{2}\left( 0,Z\right) $ such that 
\begin{equation}
\left. 
\begin{array}{c}
\widehat{J}_{\lambda ,q_{1},r_{1}}\left( h_{1},h_{2}\right) \left( k\right) =%
\left[ \widetilde{J}_{\lambda ,q_{1},r_{1}},\left( h_{1},h_{2}\right) \right]
\left( k\right) , \\ 
\forall \left( h_{1},h_{2}\right) \in H_{0}^{2}\left( 0,Z\right) \times
H_{0}^{2}\left( 0,Z\right) ,\text{ }\forall k\in \left[ k_{\min },k_{\max }%
\right] .%
\end{array}%
\right.  \label{6.12}
\end{equation}%
It follows from (\ref{6.1}), (\ref{6.3})-(\ref{6.6}) and (\ref{6.9})-(\ref%
{6.12}) that for all $k\in \left[ k_{\min },k_{\max }\right] $ 
\begin{equation*}
\left. 
\begin{array}{c}
\lim_{\left\Vert \left( h_{1},h_{2}\right) \right\Vert _{H_{0}^{2}\left(
0,Z\right) \times H_{0}^{2}\left( 0,Z\right) }}\left\Vert \left(
h_{1},h_{2}\right) \right\Vert _{H_{0}^{2}\left( 0,Z\right) \times
H_{0}^{2}\left( 0,Z\right) }^{-1}\times \\ 
\times \left\{ 
\begin{array}{c}
J_{\lambda }\left( q_{1}+h_{1},r_{1}+h_{2}\right) \left( k\right)
-J_{\lambda }\left( q_{1},r_{1}\right) \left( k\right) - \\ 
-\left[ \widetilde{J}_{\lambda ,q_{1},r_{1}},\left( h_{1},h_{2}\right) %
\right] \left( k\right)%
\end{array}%
\right\} =0.%
\end{array}%
\right.
\end{equation*}%
Hence, $\widetilde{J}_{\lambda ,q_{1},r_{1}}$ is the Fr\'{e}chet derivative
of the functional $J_{\lambda }$ at the point $\left( q_{1},r_{1}\right) ,$
i.e. 
\begin{equation}
\widetilde{J}_{\lambda ,q_{1},r_{1}}=J_{\lambda }^{\prime }\left(
q_{1},r_{1}\right) \left( k\right) \in H_{0}^{2}\left( 0,Z\right) \times
H_{0}^{2}\left( 0,Z\right) ,\text{ }\forall k\in \left[ k_{\min },k_{\max }%
\right] .  \label{6.13}
\end{equation}%
We omit the proof of the Lipschitz continuity property (\ref{5.2}) of $%
J_{\lambda }^{\prime }\left( q_{1},r_{1}\right) \left( k\right) $ since this
proof is similar with the proof of Theorem 5.3.1 of \cite{KL}.

We now prove the strong convexity property (\ref{5.3}). To do this, we use
Carleman estimate (\ref{4.3}) of Theorem 4.1. Using (\ref{4.3}), (\ref{6.7}%
)-(\ref{6.9}) and (\ref{6.13}), we obtain%
\begin{equation*}
J_{\lambda }\left( q_{1}+h_{1},r_{1}+h_{2}\right) \left( k\right)
-J_{\lambda }\left( q_{1},r_{1}\right) \left( k\right) -\left[ J_{\lambda
}^{\prime }\left( q_{1},r_{1}\right) \left( k\right) ,\left(
h_{1},h_{2}\right) \right] \geq
\end{equation*}%
\begin{equation*}
\geq \frac{1}{2}\int\limits_{0}^{Z}\left[ \left( \partial
_{z}^{2}h_{1}\right) ^{2}+\left( \partial _{z}^{2}h_{2}\right) ^{2}\right]
\varphi _{\lambda }dz-C_{1}\int\limits_{0}^{Z}\left[ \left( \partial
_{z}h_{1}\right) ^{2}+\left( \partial _{z}h_{2}\right) ^{2}\right] \varphi
_{\lambda }dz\geq
\end{equation*}%
\begin{equation}
\geq \frac{1}{2}C_{0}\int\limits_{0}^{Z}\left[ \left( \partial
_{z}^{2}h_{1}\right) ^{2}+\left( \partial _{z}^{2}h_{2}\right) ^{2}\right]
\varphi _{\lambda }dz+  \label{6.14}
\end{equation}%
\begin{equation*}
+\frac{1}{2}C_{0}\lambda \int\limits_{0}^{Z}\left[ \left( \partial
_{z}h_{1}\right) ^{2}+\left( \partial _{z}h_{2}\right) ^{2}+\lambda
^{2}\left( h_{1}^{2}+h_{2}^{2}\right) \right] {\varphi }_{\lambda }\emph{dz}-
\end{equation*}%
\begin{equation*}
-C_{1}\int\limits_{0}^{Z}\left[ \left( \partial _{z}h_{1}\right) ^{2}+\left(
\partial _{z}h_{2}\right) ^{2}\right] \varphi _{\lambda }dz,\text{ }\forall
\lambda \geq \lambda _{0}.
\end{equation*}%
Hence, we can choose a sufficiently large number $\lambda _{1}=\lambda
_{1}\left( R,Z,k_{\min },k_{\max },\varepsilon \right) \geq \lambda _{0}$
such that (\ref{6.14}) becomes%
\begin{equation*}
J_{\lambda }\left( q_{1}+h_{1},r_{1}+h_{2}\right) \left( k\right)
-J_{\lambda }\left( q_{1},r_{1}\right) \left( k\right) -\left[ J_{\lambda
}^{\prime }\left( q_{1},r_{1}\right) \left( k\right) ,\left(
h_{1},h_{2}\right) \right] \geq
\end{equation*}%
\begin{equation*}
\geq \frac{1}{2}C_{0}\int\limits_{0}^{Z}\left[ \left( \partial
_{z}^{2}h_{1}\right) ^{2}+\left( \partial _{z}^{2}h_{2}\right) ^{2}\right]
\varphi _{\lambda }dz+
\end{equation*}%
\begin{equation*}
+C_{1}\lambda \int\limits_{0}^{Z}\left[ \left( \partial _{z}h_{1}\right)
^{2}+\left( \partial _{z}h_{2}\right) ^{2}+\lambda ^{2}\left(
h_{1}^{2}+h_{2}^{2}\right) \right] \varphi _{\lambda }dz\geq
\end{equation*}%
\begin{equation*}
\geq C_{1}e^{-2\lambda Z}\left\Vert \left( h_{1},h_{2}\right) \right\Vert
_{H^{2}\left( 0,Z\right) \times H^{2}\left( 0,Z\right) }^{2}=
\end{equation*}%
\begin{equation*}
=C_{1}e^{-2\lambda Z}\left\Vert \left( q_{2},r_{2}\right) \left( k\right)
-\left( q_{1},r_{1}\right) \left( k\right) \right\Vert _{H^{2}\left(
0,Z\right) \times H^{2}\left( 0,Z\right) }^{2},\text{ }\forall \lambda \geq
\lambda _{1},
\end{equation*}%
which proves (\ref{5.3}).

Existence and uniqueness of the minimizer $\left( q_{\min ,\lambda }\left(
z,k\right) ,r_{\min ,\lambda }\left( z,k\right) \right) \in \overline{%
B\left( R\right) }$\ of the functional $J_{\lambda }\left( q,r\right) \left(
k\right) $\ on the set $\overline{B\left( R\right) }$ as well as inequality (%
\ref{5.4}) easily follow immediately from (\ref{5.3}) and a combination of
Lemma 5.2.1 with Theorem 5.2.1 of \cite{KL}. $\square $

\section{Proof of Theorem 5.2}

\label{sec:7}

Let $I_{\lambda }\left( q,r\right) $ be the functional defined in (\ref{5.12}%
). As stated in lines below (\ref{5.12}), an obvious analog of Theorem 5.1
is valid for $I_{\lambda }\left( q,r\right) $ for values of the parameter $%
\lambda $ as in (\ref{5.13}). Recall that $\left( \widetilde{q}_{\min
,\lambda },\widetilde{r}_{\min ,\lambda }\right) \in \overline{B_{0}\left(
2R\right) }$\emph{\ }is the unique minimizer on the set $\overline{%
B_{0}\left( 2R\right) }$ of the functional $I_{\lambda }\left( q,r\right) $
for $\lambda \geq \lambda _{2}.$ By (\ref{5.3}), (\ref{5.10}) and the second
line of (\ref{5.12})%
\begin{equation*}
I_{\lambda }\left( \widetilde{q}^{\ast },\widetilde{r}^{\ast }\right) \left(
k\right) -I_{\lambda }\left( \widetilde{q}_{\min ,\lambda },\widetilde{r}%
_{\min ,\lambda }\right) \left( k\right) -
\end{equation*}%
\begin{equation}
-\left[ I_{\lambda }^{\prime }\left( \widetilde{q}_{\min ,\lambda },%
\widetilde{r}_{\min ,\lambda }\right) \left( k\right) ,\left( \widetilde{q}%
^{\ast }-\widetilde{q}_{\min ,\lambda },\widetilde{r}^{\ast }-\widetilde{r}%
_{\min ,\lambda }\right) \left( k\right) \right] \geq  \label{7.1}
\end{equation}%
\begin{equation*}
\geq C_{1}e^{-2\lambda Z}\left\Vert \left( \widetilde{q}^{\ast },\widetilde{r%
}^{\ast }\right) \left( k\right) -\left( \widetilde{q}_{\min ,\lambda },%
\widetilde{r}_{\min ,\lambda }\right) \left( k\right) \right\Vert
_{H^{2}\left( 0,Z\right) \times H^{2}\left( 0,Z\right) }^{2},\text{ }\forall
k\in \left[ k_{\min },k_{\max }\right] .
\end{equation*}%
By (\ref{5.4})%
\begin{equation*}
-\left[ I_{\lambda }^{\prime }\left( \widetilde{q}_{\min ,\lambda },%
\widetilde{r}_{\min ,\lambda }\right) \left( k\right) ,\left( \widetilde{q}%
^{\ast }-\widetilde{q}_{\min ,\lambda },\widetilde{r}^{\ast }-\widetilde{r}%
_{\min ,\lambda }\right) \left( k\right) \right] \leq 0.
\end{equation*}%
In addition, obviously $-I_{\lambda }\left( \widetilde{q}_{\min ,\lambda },%
\widetilde{r}_{\min ,\lambda }\right) \left( k\right) \leq 0.$ Hence, (\ref%
{7.1}) implies%
\begin{equation*}
I_{\lambda }\left( \widetilde{q}^{\ast },\widetilde{r}^{\ast }\right) \left(
k\right) \geq
\end{equation*}%
\begin{equation}
\geq C_{1}e^{-2\lambda Z}\left\Vert \left( \widetilde{q}^{\ast },\widetilde{r%
}^{\ast }\right) \left( k\right) -\left( \widetilde{q}_{\min ,\lambda },%
\widetilde{r}_{\min ,\lambda }\right) \left( k\right) \right\Vert
_{H^{2}\left( 0,Z\right) \times H^{2}\left( 0,Z\right) }^{2},\text{ }\forall
k\in \left[ k_{\min },k_{\max }\right] .  \label{7.2}
\end{equation}%
Next, by (\ref{5.10}) and (\ref{5.12}) 
\begin{equation*}
I_{\lambda }\left( \widetilde{q}^{\ast },\widetilde{r}^{\ast }\right) \left(
k\right) =J_{\lambda }\left( \widetilde{q}^{\ast }+F_{1},\widetilde{r}^{\ast
}+F_{2}\right) \left( k\right) =
\end{equation*}%
\begin{equation}
=J_{\lambda }\left( \left( \widetilde{q}^{\ast }+F_{1}^{\ast }\right)
+\left( F_{1}-F_{1}^{\ast }\right) ,\left( \widetilde{r}^{\ast }+F_{2}^{\ast
}\right) +\left( F_{2}-F_{2}^{\ast }\right) \right) \left( k\right) =
\label{7.3}
\end{equation}%
\begin{equation*}
=J_{\lambda }\left( q^{\ast }+\left( F_{1}-F_{1}^{\ast }\right) ,r^{\ast
}+\left( F_{2}-F_{2}^{\ast }\right) \right) .
\end{equation*}%
Now, $J_{\lambda }\left( q^{\ast },r^{\ast }\right) =0.$ Hence, using (\ref%
{5.5}) and (\ref{7.3}), we obtain%
\begin{equation*}
I_{\lambda }\left( \widetilde{q}^{\ast },\widetilde{r}^{\ast }\right) \left(
k\right) =J_{\lambda }\left( q^{\ast }+\left( F_{1}-F_{1}^{\ast }\right)
,r^{\ast }+\left( F_{2}-F_{2}^{\ast }\right) \right) \leq C_{1}\delta ^{2}.
\end{equation*}%
Combining this with (\ref{7.2}) and setting $\lambda =\lambda _{2},$ we
obtain%
\begin{equation}
\left\Vert \left( \widetilde{q}^{\ast },\widetilde{r}^{\ast }\right) \left(
k\right) -\left( \widetilde{q}_{\min ,\lambda _{2}},\widetilde{r}_{\min
,\lambda _{2}}\right) \left( k\right) \right\Vert _{H^{2}\left( 0,Z\right)
\times H^{2}\left( 0,Z\right) }\leq C_{1}\delta ,\text{ }\forall k\in \left[
k_{\min },k_{\max }\right] .  \label{7.4}
\end{equation}%
Using (\ref{5.10}) and (\ref{5.14}), we obtain%
\begin{equation*}
\widetilde{q}^{\ast }-\widetilde{q}_{\min ,\lambda _{2}}=\left( \widetilde{q}%
^{\ast }+F_{1}^{\ast }\right) -\left( \widetilde{q}_{\min ,\lambda
_{2}}+F_{1}\right) -\left( F_{1}^{\ast }-F_{1}\right) =
\end{equation*}%
\begin{equation}
=\left( q^{\ast }-\overline{q}_{\min ,\lambda _{2}}\right) -\left(
F_{1}^{\ast }-F_{1}\right) .  \label{7.5}
\end{equation}%
Similarly 
\begin{equation}
\widetilde{r}^{\ast }-\widetilde{r}_{\min ,\lambda _{2}}=\left( r^{\ast }-%
\overline{r}_{\min ,\lambda _{2}}\right) -\left( F_{2}^{\ast }-F_{2}\right) .
\label{7.6}
\end{equation}%
Hence, using (\ref{5.5}), (\ref{7.4})-(\ref{7.6}) and triangle inequality,
we obtain%
\begin{equation}
\left\Vert \left( q^{\ast },r^{\ast }\right) \left( k\right) -\left( 
\overline{q}_{\min ,\lambda _{2}},\overline{r}_{\min ,\lambda _{2}}\right)
\left( k\right) \right\Vert _{H^{2}\left( 0,Z\right) \times H^{2}\left(
0,Z\right) }\leq C_{1}\delta ,\text{ }\forall k\in \left[ k_{\min },k_{\max }%
\right] .  \label{7.7}
\end{equation}%
\begin{equation}
I_{\lambda }\left( \widetilde{q}_{\min ,\lambda },\widetilde{r}_{\min
,\lambda }\right) =\min_{\overline{B_{0}\left( 2R\right) }}I_{\lambda
}\left( q,r\right) .5.130
\end{equation}%
\emph{\ }Hence, using (\ref{5.5}) and (\ref{7.7}), we obtain (\ref{5.140}).
By (\ref{5.9}), (\ref{5.100}), (\ref{5.12}) and (\ref{5.130})%
\begin{equation*}
I_{\lambda _{2}}\left( \widetilde{q}_{\min ,\lambda _{2}},\widetilde{r}%
_{\min ,\lambda _{2}}\right) =J_{\lambda _{2}}\left( \widetilde{q}_{\min
,\lambda _{2}}+F_{1},\widetilde{r}_{\min ,\lambda _{2}}+F_{2}\right) =
\end{equation*}%
\begin{equation}
=J_{\lambda _{2}}\left( \overline{q}_{\min ,\lambda _{2}},\overline{r}_{\min
,\lambda _{2}}\right) \leq J_{\lambda _{2}}\left( q,r\right) =  \label{7.8}
\end{equation}%
\begin{equation*}
=J_{\lambda _{2}}\left( \left( q-F_{1}\right) +F_{1},\left( r-F_{2}\right)
+F_{2}\right) ,\text{ }\forall \left( q,r\right) \in \overline{B\left(
R\right) }.
\end{equation*}%
Let $\left( q_{\min ,\lambda _{2}},r_{\min ,\lambda _{2}}\right) \in 
\overline{B\left( R\right) }$ be the unique minimizer of the functional $%
J_{\lambda _{2}}\left( q,r\right) $ on the set $\overline{B\left( R\right) }%
, $ the existence of which is guaranteed by Theorem 5.1. Hence, (\ref{7.8})
implies that 
\begin{equation}
J_{\lambda _{2}}\left( \overline{q}_{\min ,\lambda _{2}},\overline{r}_{\min
,\lambda _{2}}\right) \leq J_{\lambda _{2}}\left( q_{\min ,\lambda
_{2}},r_{\min ,\lambda _{2}}\right) .  \label{7.9}
\end{equation}%
However, since by (\ref{5.140}) $\left( \overline{q}_{\min ,\lambda _{2}},%
\overline{r}_{\min ,\lambda _{2}}\right) \in \overline{B\left( R\right) },$
then (\ref{7.7}) and (\ref{7.9}) imply (\ref{5.141}) and (\ref{5.15}) and (%
\ref{5.150}). $\square $

\begin{center}
\textbf{Acknowledgment}
\end{center}

This work was partially supported by US National Science Foundation grant
DMS 2436227.

\end{document}